\documentclass[11pt]{amsart}

\input epsf
\usepackage{latexsym}
\usepackage{amssymb,amsmath,amsthm,amscd, amssymb,
graphicx, enumerate, verbatim}
\usepackage[all]{xy}
\numberwithin{equation}{section}
\newtheorem{lem}[equation]{Lemma}
\newtheorem{prop}[equation]{Proposition}

\newtheorem{Example}[equation]{Example}

\newtheorem{remark}[equation]{Remark}
\newenvironment{rmk}{\begin{remark}\rm}{\end{remark}}
\def\co{\colon\thinspace}

\newcommand{\Ric}{\mbox{Ric}}

\def\parT#1{\frac{\partial #1}{\partial t}}

\begin{document}
\abovedisplayskip=6pt plus3pt minus3pt \belowdisplayskip=6pt
plus3pt minus3pt
\title {Curvature bounds via Ricci smoothing }
\rm
\thanks{\it This work was supported in part by the NSF grant \# DMS-0204187}
\author{Vitali Kapovitch}
\address{Vitali Kapovitch\\Department of Mathematics\\University of California\\
Santa Barbara, CA 93110}\email{vitali@math.ucsb.edu}
\date{}
\begin{abstract}
We give a proof of the fact that  the upper and the lower
sectional curvature bounds of a complete manifold vary at a
bounded rate under the Ricci flow.
\end{abstract}
\maketitle
Let $(M^n,g)$ be a complete Riemannian manifold with $|\sec
(M)|\le 1$.
 Consider the Ricci flow of $g$ given by
 \begin{equation}\label{ricflow}
\frac{\partial}{\partial t}g=-2\Ric (g)
 \end{equation}
 It is  known ( see \cite{Ha},~\cite{Shi}) that (\ref{ricflow}) has a solution on  $[0, T]$
 for some $T>0$.
It is  also known (see  \cite{BM,Shi}) that the
solution smoothes out the metric. Namely, $g_t$   satisfies
\begin{equation}\label{e1:smooth}
e^{-c(n)t}g\le g_t\le e^{c(n)t}g\quad |\nabla-\nabla_t|\le
c(n)t\quad |\nabla^m R_{ijkl}(t)|\le c(n,m,t)
\end{equation}
Moreover, by \cite{Shi} , sectional curvature of $g(t)$
satisfies
\begin{equation}\label{est-sec1}
|K_{g_t}|\le C(n,T)
\end{equation}
 This result proved to
be a very useful technical tool in many situations and in
particular in the theory of convergence with two-sided curvature
bounds ( see \cite{CFG,Rong,PT} etc). However it turns
out that in applications to convergence with two-sided curvature
bounds  in addition to the
above properties it is often convenient to know that  $\sup K_{g_t}$ and $\inf K_{g_t}$
also vary at the bounded rate and in particular, the upper and the
lower curvature bounds for $g_t$ are almost the same as for $g$
for sufficiently small $t$.
For example, it is very useful to know that if $g_{0}$ has pinched
positive~\cite{Rong} or negative~\cite{Kan,BK} curvature, then
$g_{t}$ has almost the same pinching.

This fact has apparently been known to
some experts and it was used without a proof by various people
(see e.g \cite{Kan}).
 A careful proof was given in \cite{Rong} in case of a compact $M$.
 To the best of our knowledge, no proof exists in the literature in case of a
 noncompact $M$.
 The purpose of this note is to rectify this situation.
  To this end we prove
  \begin{prop}\label{prop:sec}
In the above situation one has
\[
\inf K_{g}-C(n,T)t\le K_{g_t}\le \sup K_{g}+C(n,T)t
\]
  \end{prop}
  \begin{proof}
Throughout the proof we will denote by $C$ various constants depending only on $n,T$.
The proof in \cite{Rong} relies on the maximum principle applied
to the evolution equation for the curvature tensor $Rm$ which  can
be computed to have the form \cite{Shi}
\begin{equation}\label{ev}
\frac{\partial}{\partial t}R_{ijkl}=\Delta R_{ijkl} + P(Rm)
\end{equation}
where $P(Rm)$ is a quadratic polynomial in $Rm$.
However, in noncompact case the maximum principle can not be
applied directly. We will use a local version of the maximum
principle often employed in \cite{Shi}.
Let $\chi\co \mathbb R\to \mathbb R$ be a smooth function
satisfying
\begin{enumerate}[(1)]
\item $\chi\ge 0$ and is nonincreasing
\item $\chi (x)=
\begin{cases}
1 \textrm{ for } x\le 1\\
\textrm{nonincreasing for } 1\le x\le 2\\
0\textrm{ for } x\ge 2
\end{cases}$
\item $|\chi''(x)|\le 8$
\item $\left|\frac{(\chi'(x))^2}{\chi(x)}\right|\le 16$
\end{enumerate}
Fix $z\in M$ and let $d_z(x,t)=d_{g_t}(x, z)$ be the distance with
respect to $g_t$. Put $\xi_z(x,t)=\chi(d_z(x,t))$ . Using the
properties of $\chi$  we obtain
\begin{enumerate}[(i)]
\item $0\le \xi_z\le 1$
\item $|\nabla\xi_z|\le C$
\item \label{est:xi}$\Delta \xi_z\ge C$ in the barrier sense
\item $\frac{|\nabla\xi_z|^2}{|\xi_z|}\le C$
\item $|\frac{\partial\xi_z(x,t)}{\partial t}|\le C$.
\end{enumerate}
To see (\ref{est:xi}) we compute $\Delta \xi_z=\chi''(d_z)|\nabla d_z|^2+\chi'(d_z)\Delta d_z\ge C$ because$\chi'\le 0$
and $\Delta d_z\le C$ for $d_z\ge 1$ by Laplace comparison for spaces with  $\sec\ge -1$.
Finally, (v) holds by the evolution equation of the metric (\ref{ricflow}) and the estimate
(\ref{est-sec1}).

Assume for now that $\sup K_{g_{t}}\ge 0$ for all $t\in [0,T]$.
Let $\bar{A}(t)=\sup K_{g_{t}}$ and  $\bar{A}_z(t)=\max_{(x,\sigma)}\{K_{g_t}(x,\sigma)\xi_z(x,t)\}$
where $x\in M$, $\sigma$ is a 2-plane at $x$. Clearly $\bar{A}(t)=\sup_z \bar{A}_z(t)$.

We want to show that $\bar{A}'_z(t)\le C$ independent of $z,t$. Fix $t_0\in [0,T]$ and
 let $\phi_z(x,\sigma,t)=K_{g_{t}}(x,\sigma)\xi_z(x,t)$.
By a standard argument,  it is enough to check that
$\parT{\phi_z}(x_0,\sigma_0,t_0)\le C$ for any point of maximum of
$\phi_z(\cdot, t_0)$.

Let $U,V$ be a basis of $\sigma_0$ orthonormal with respect to $g_{t_0}$.
Extend $U,V$ to  constant vector fields in normal coordinates at $x_0$ with respect to $g_{t_0}$.

Let
$\Phi_z(x,t)=K_{g_t}(x,U,V)\xi_z(x)=\frac{Rm(t)(U,V,U,V)}{|U\wedge
V|^2_{g_t}}\xi_z(x)$.

It is easy to see (cf.~\cite{Rong}) that
\begin{equation}\label{e:UV}
|U\wedge V(x_0)|_{g_t}\le C,|\nabla |U\wedge V(x_0)|_{g_t}|\le C \textrm{ and }
|\nabla^2 |U\wedge V(x_0)|_{g_t}|\le C
\end{equation}
By construction, $\Phi_z(x,t_0)$ has a local maximum at $x_0$ and
$\parT{ \phi_z(x_0,\sigma_0,t_0)}=\parT{ \Phi_z(x_0,t_0)}$.
Therefore $\nabla \Phi_z(x_0,t_0)=0$ and
 $\Delta \Phi_z(x_0,t_0)\le 0$.
We compute
\begin{equation}\label{e:ev3}
\begin{split}
\parT{ \Phi_z(x_0,t_0)}=\Delta \Phi_z(x_0,t_0)- Rm(x_0,t_0)(U,V,U,V)\xi_z(x_0,t_0)\parT{}\left({\frac{1}{|U\wedge V|^2}}
\right)\\
-2\nabla {Rm(x_0,t_0)(U,V,U,V)}\nabla\left(\frac{\xi_z(x_0,t_0)}{|U\wedge V|^2}\right)
-Rm(x_0,t_0)(U,V,U,V) \Delta \left(\frac{\xi_z(x_0,t_0)}{|U\wedge V|^2}\right)-\\
\frac{P(Rm(x_0,t_0))\xi_z(x_0,t_0)}{|U\wedge V|^2}-K_{g_t}(x,U,V)\parT{\xi_z(x_0,t_0)}
\end{split}
\end{equation}
We claim that the RHS is bounded above by $C$. The only terms that need explaining are the third and the forth
summands. Let $f(x)=\frac{\xi_z(x,t_0)}{|U\wedge V|^2}$.

 To see that the third term is bounded we observe that $\nabla \Phi_z(x_0,t_0)=0$ yields
$\nabla Rm(x_0,t_0)(U,V,U,V)f(x_0)+Rm(x_0,t_0)(U,V,U,V)\nabla f(x_0)=0$,\\ $
\nabla Rm(x_0,t_0)(U,V,U,V)=-\frac{\nabla f(x_0)}{f(x_0)}Rm(x_0,t_0)(U,V,U,V)$ and hence \\
$|\nabla Rm(x_0,t_0)(U,V,U,V)\nabla f(x_0)|\le C$ by the property (iv) of $\xi_z$ above.
The fourth term is bounded because $\Delta f=\Delta \xi_z(x_0) \frac{1}{|U\wedge V|^2}+
2\nabla \xi_z(x_0) \nabla\left(\frac{1}{|U\wedge V|^2}\right)+\xi_z(x_0) \Delta\left(\frac{1}{|U\wedge V|^2}\right)\ge C$
by (\ref{e:UV}) and the property (iii) of $\xi_z$.
 Thus by (\ref{e:ev3}) we have $\parT{\phi_z}(x_0,\sigma_0,t_0)=\parT{ \Phi_z(x_0,t_0)}\le C$.
Thus $\bar{A}'_z(t)\le C$ for all $z\in M, t\in[0,T]$ and hence
$\bar{A}'(t)\le C$ for all $t\in[0,T]$ This concludes the proof in
the case $\sup K_{g_{t}}\ge 0$. The general case can be easily
reduced to this one by replacing the function
$K_{g_{t_0}}(x,\sigma)$ by $K_{g_{t_0}}(x,\sigma)+C$. The argument
for $\inf K_{g_{t}}$ is the same except there we  can actually
always assume that $\inf K_{g_{t}}\le 0$ since otherwise the
manifold $M$ is compact and our statement is known by
~\cite{Rong}.
\end{proof}
\begin{rmk}\label{r:1}
 By changing the cutoff function $\xi_z(\cdot )$ to $\chi(
d(\cdot,z)/R)$ in the proof of Proposition~\ref{prop:sec} we see
that the same proof actually shows that the {\it local} maximum
and minimum of the curvature vary linearly. Namely, under
condition of the Proposition, for any $R>0$ there exists
$C=C(T,R)$ such that for any $z\in M$ we have
\[
\inf_{B(z,R)} K_{g}-C(n,R,T)t\le K_{g_t}|_{B(z,R)}\le \sup_{B(z,R)} K_{g}+C(n,R,T)t
\]
However, as constructed, $C(n,R,T)\to\infty$ as $R\to 0$.
\end{rmk}
\begin{rmk}
A slightly more careful examination of the proof of
Proposition~\ref{prop:sec} shows that the local rate of change of
the curvature bounds is proportional to the local absolute
curvature bounds, i.e $\bar{A}'_z(t)\le C(n,T)\cdot sup_{x\in
B(z,2)}|Rm(x)|$. In particular,  if $(M^n, g)$ is asymptotically
flat then so is $(M^n, g_t)$ and it has the same curvature decay
rate as $(M^n, g)$.
\end{rmk}
\bibliographystyle{amsalpha}
\def\cprime{$'$} \def\cprime{$'$}
\providecommand{\bysame}{\leavevmode\hbox to3em{\hrulefill}\thinspace}
\providecommand{\MR}{\relax\ifhmode\unskip\space\fi MR }
\providecommand{\MRhref}[2]{%
  \href{http://www.ams.org/mathscinet-getitem?mr=#1}{#2}
}
\providecommand{\href}[2]{#2}

\end{document}